\documentclass[11pt,reqno]{amsart}
\usepackage{amsfonts,amsmath,amssymb}
\newtheorem{prop}{Proposition}

\DeclareMathOperator{\arcosh}{arcosh}

\def\beq#1#2\eeq{%
        \begin{equation}%
        \label{#1}%
            #2%
        \end{equation}%
    }

\usepackage{color}
\usepackage{graphicx}
\usepackage{epstopdf}

\theoremstyle{plain}
\newtheorem{theorem}{Theorem}

\theoremstyle{remark}

\theoremstyle{definition}

\def\Z{\mathbb{Z}}
\def\N{\mathbb{N}}
\def\R   {\mathbb{R}}
\def\T   {\mathbb{T}}
\def\Q{\mathbb{Q}}

\def \cT {\mathcal T}

\def\g{\gamma}

\def\eps{\varepsilon}

\title[Fock and Mather]{Markov numbers, Mather's $\beta$ function and stable norm}

\author{A. Sorrentino}\address{Dipartimento di Matematica, Universit\`a degli Studi di Roma Tor Vergata, Rome,
Italy}
\email{sorrentino@mat.uniroma2.it}

\author{A.P. Veselov}
\address{Department of Mathematical Sciences,
Loughborough University, Loughborough LE11 3TU, UK  and Moscow State University, Moscow 119899, Russia}
\email{A.P.Veselov@lboro.ac.uk}

\begin{document}

\maketitle

\begin{abstract}
V. Fock \cite{Fock} introduced an interesting function $\psi(x), x \in \mathbb R$ related to Markov numbers.
We explain its relation to {Federer-Gromov's} stable norm and Mather's $\beta$-function, and use this to study its properties. 
We prove that $\psi$ and its natural generalisations are differentiable at every irrational $x$ and non-differentiable otherwise, {by exploiting} the relation with {length of} simple closed geodesics on the punctured {or} one-holed tori with the hyperbolic metric and the results by Bangert {\cite{Bangert2}} and McShane-Rivin {\cite{McShane}}.
\end{abstract}

\section{Introduction}

In 1880 A.A. Markov \cite{M80} 
discovered a remarkable relation between the theory of binary quadratic forms and the following Diophantine equation known as {\it Markov equation}
\begin{equation}
\label{Markov}
x^2 + y^2+ z^2 = 3xyz.
\end{equation}
Markov showed that all positive integer solutions (known as {\it Markov triples}) can be {obtained} from the obvious $(1,1,1)$ by applying the {\it Vieta symmetry}
\begin{equation}
\label{Markovinv}
(x, y, z) \; {\longmapsto} \; (x, y, 3xy-z)
\end{equation}
 and permutations. The elements of Markov triples are the famous {\it Markov numbers}
\[1, 2, 5, 13, 29, 34, 89, 169, 194, 233, 433, 610, 985,...\]
which play a very important role in number theory \cite{Aigner}, in the theory of Frobenius manifolds and related Painlev\'e-VI equation, Teichm\"uller spaces and algebraic geometry (see \cite{SV} and references therein).

The Markov numbers can be naturally labelled by  {rationals} $x \in [0,1/2]$ using the Farey tree  (see fig. 1). Recall that at each vertex of the Farey tree we have fractions $\frac{a}{b}$, $\frac{c}{d}$ and their {\it Farey mediant} $\frac{a+c}{b+d}$ (see e.g. \cite{SV}). On the Markov tree the triples at two neighbouring vertices are related by the involution (\ref{Markovinv}).

Let $m\left(\frac{p}{q} \right)$ {be} the Markov number corresponding to $\frac{p}{q}$ on {the} Farey tree. The function $m(\frac{p}{q})$ can be extended to all {rationals} $\frac{p}{q}$ using the symmetry
$$
m(1-x)=m(1/x)=m(x).
$$

\begin{figure}[h]
\label{fig:CFTree}
\begin{center}
\includegraphics[trim = 0mm 30mm 0mm 22mm, clip, height=38mm]{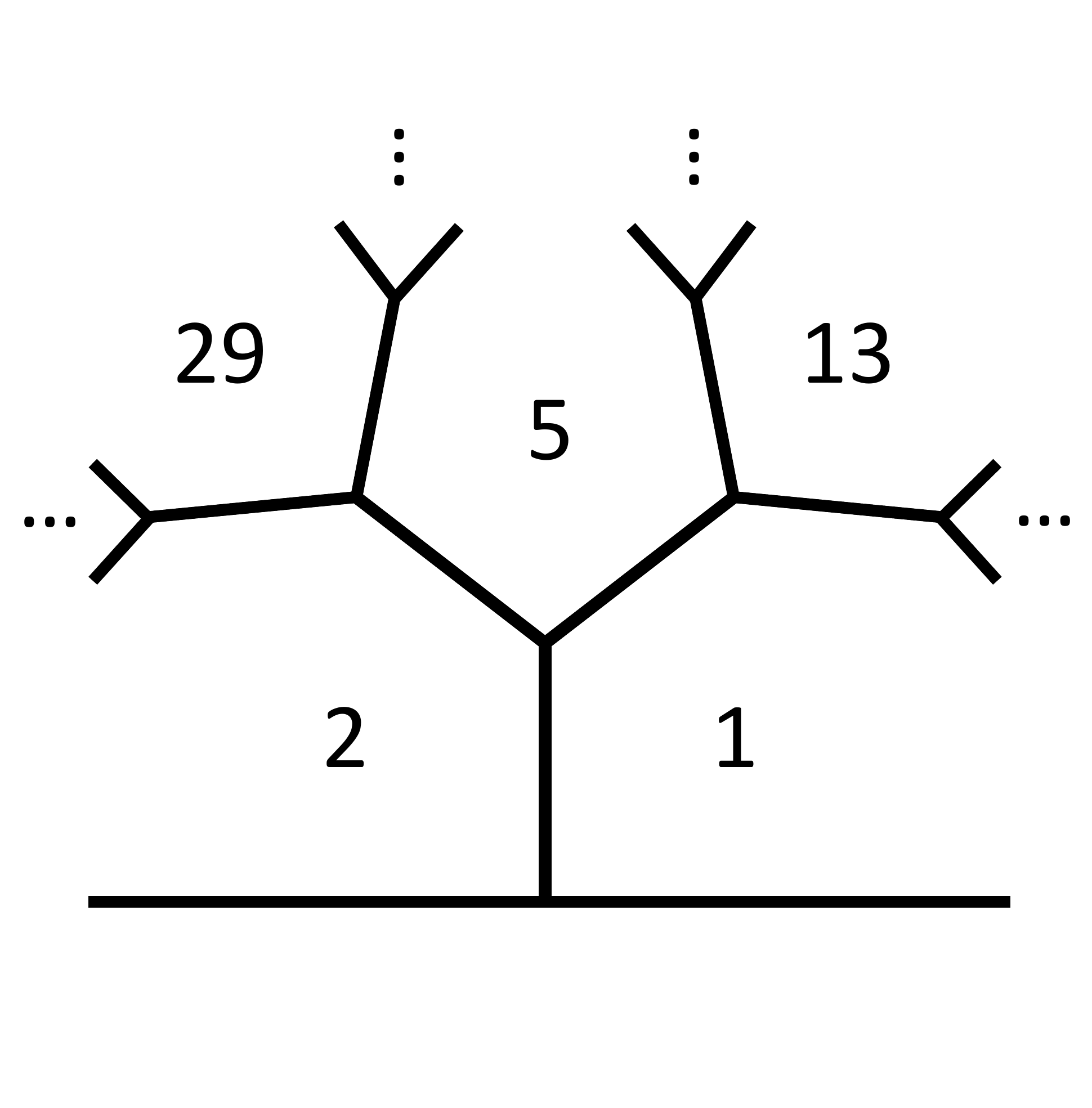}  \hspace{8pt} \quad  \includegraphics[trim = 0mm 30mm 0mm 30mm, clip, height=38mm]{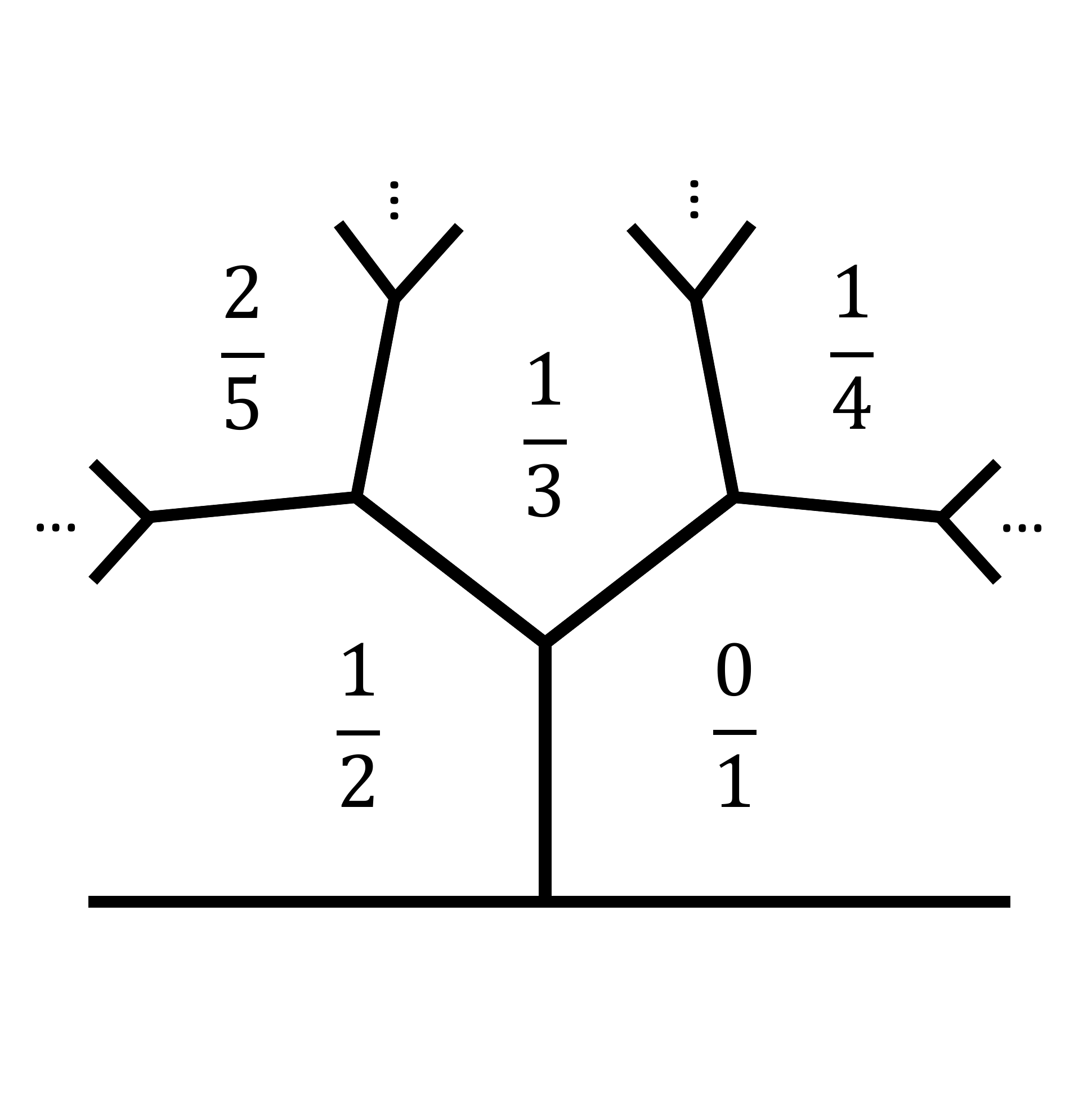}
\caption{\small Markov tree and corresponding branch of Farey tree}
\end{center}
\end{figure}

Following V. Fock \cite{Fock} consider the following function $\psi(x), \, x \in \mathbb R.$ {At a} rational $x = \frac{p}{q}$ {this function} is defined as
\begin{equation}\label{fock}
\psi \left(\frac{p}{q}\right) {:=} \frac{1}{q} \arcosh \left(\frac{3}{2} m \left(\frac{p}{q}\right)\right).\\
\end{equation}

\medskip

\begin{theorem} [Fock]
Function $\psi$ can be extended to a continuous convex function {on the whole} $\mathbb R.$ 
\end{theorem}

The aim of this note is to explain the relation of this result with {{Aubry-Mather theory} (more specifically with the so-called Mather's $\beta$-function \cite{Ma90,Ma, Sorrentinobook}} and {Federer-Gromov's} stable norm \cite{GLP})\footnote{{
This concept appeared for the first time in Federer \cite{Federer} and was named {\it stable norm} in Gromov \cite{GLP}}}.
We use the well-known relation {between} Markov numbers {and} the lengths of simple closed geodesic on a punctured torus (see \cite{Cohn, Series}).
This allows us to prove the following \\

\begin{theorem}
\label{main}
Fock's function $\psi$ is differentiable at all {irrationals} and not differentiable at {any rational}.\\
\end{theorem}

We also show also that these results can be naturally generalised {to} the solutions of the Diophantine equation \cite{SV}
$$
X^2+Y^2+Z^2=XYZ+4-4a^6, \qquad a \in \mathbb N,
$$ 
{which are} related to  hyperbolic tori with a hole (see  next section). 

%\medskip

\section{Markov equation and hyperbolic tori}

We explain now in more detail the relation of Markov numbers with the closed geodesics on the punctured torus with hyperbolic metric, which was found by  Gorshkov \cite{Gorshkov} in his thesis in 1953 and, independently, by Cohn \cite{Cohn} (see also \cite{Haas, Series}).

The fundamental group of the punctured torus is the free group $F_2$. The hyperbolic structure corresponds to a realisation of $F_2$ as a discrete (Fuchsian) subgroup of $SL_2(\mathbb R)$. Let $A,B \in SL_2(\mathbb R)$ be the corresponding generators.
We have the classical {\it Fricke identities}: for any $A,B \in SL_2(\mathbb R), C=AB$ 
$$
tr\, AB + tr\, AB^{-1}=tr\, A\, tr\, B,
$$
\begin{equation}
\label{fricke}
(tr\, A)^2+(tr\, B)^2+(tr\, C)^2=tr\, A \, tr\, B\, tr\,C +tr\, (ABA^{-1}B^{-1})+2.
\end{equation}

The puncture corresponds to the condition
$$
tr\, (ABA^{-1}B^{-1})=-2,
$$
so the Teichm\"uller space of punctured tori is represented by the {\it real} Markov surface
\begin{equation}
\label{xyz}	
X^2 + Y^2 + Z^2=XYZ,\quad X,Y,Z \in \mathbb R,
\end{equation}
where $X=tr\, A,\, Y=tr\, B,\, Z=tr\, C$.

The corresponding  mapping class group $GL_2(\mathbb Z)$ acts by {permutations} of $X,Y,Z$ and {by} Vieta involution
$$
(X,Y,Z)  \;{\longmapsto}\; (X,Y, XY-Z).
$$
The orbit starting from the symmetric solution $(3,3,3)$ simply consists  of 
 Markov triples multiplied by 3: $$X=3x, \,Y=3y,\, Z=3z.$$
It corresponds to the punctured equianharmonic hyperbolic torus with 3-fold symmetry, which implies that the Markov numbers are related to the lengths $l$ of simple closed geodesics by
$$m=\frac{2}{3} \cosh \frac{l}{2}.$$ 
The corresponding matrices $A,B$ can be chosen as
\begin{equation}
\label{initial}
A=\left(\begin{array}{cc} 1 & 1 \\ 1 & 2
\end{array}\right), \quad B=\left(\begin{array}{cc} 3 & 4 \\ 2 & 3
\end{array}\right)
\end{equation}
and {they} generate the commutator subgroup of $SL_2(\mathbb Z).$

Consider the following generalisation of these matrices proposed in \cite{SV} in relation with the computation of the Lyapunov exponents of Markov and Euclid trees:
\begin{equation}
\label{cohna}
A_a = \begin{pmatrix}
  1-a+a^2 & a^2 \\
  a & a+1 \\
\end{pmatrix}, \,\, B_{a} = 
\begin{pmatrix}
  1-2a+4a^2 & 4a^2 \\
  2a & 2a+1 \\
\end{pmatrix}, \quad a \in \mathbb N.
\end{equation}
The trace of the corresponding commutator is $2-4a^6$, so we have a hyperbolic torus with a hole of length $l=2\arcosh (2a^6-1).$
When $a=1$ we have the punctured torus with $l=0$.

The corresponding traces are solutions of the Diophantine equation \begin{equation}
\label{xyza}
X^2 + Y^2 + Z^2=XYZ+4-4a^6,
\end{equation}
which is a particular case of the one studied by Mordell  \cite{Mordell}. It has no fully symmetric integer solutions, but has a solution with $X=Y$:
\begin{equation}\label{orbit}
X=Y=a^2+2, \, Z=4a^2+2.
\end{equation}
Applying the permutations and Vieta involution we have the following generalisation of the Markov tree \cite{SV}.

\begin{figure}[h]
\centering
\includegraphics[trim = 0mm 30mm 0mm 30mm, clip, height=38mm]{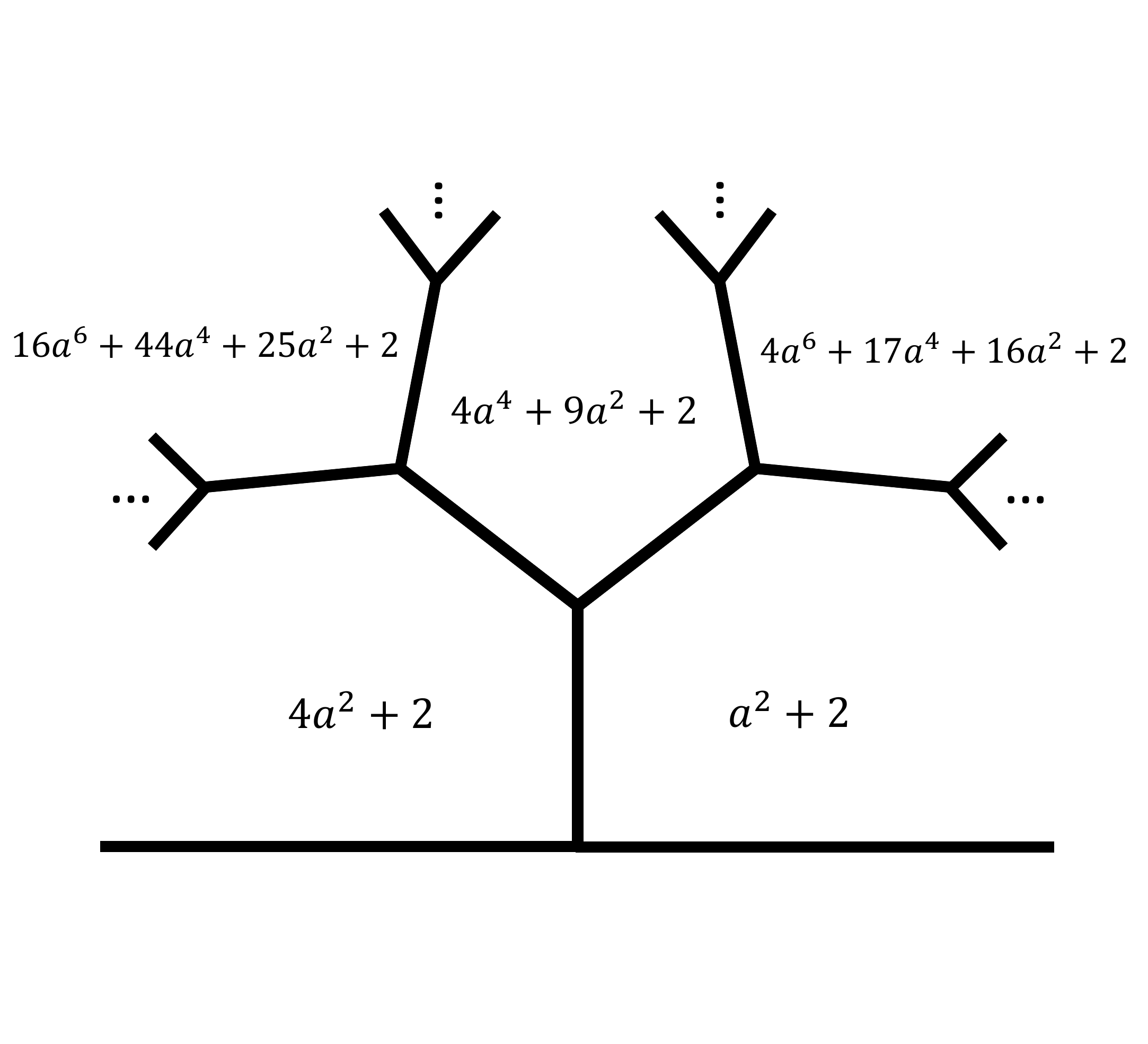} \hspace{8pt} \quad  \includegraphics[trim = 0mm 30mm 0mm 30mm, clip, height=38mm]{FareyTree21}
\label{fig:M-a-tree}
\caption{The $a$-generalisation of Markov tree and corresponding Farey fractions.}
\end{figure}

Let $X(\frac{p}{q})$, {where} $\frac{p}{q} \in [0,\frac{1}{2}]$, be the $a$-Markov number corresponding to Farey fraction $\frac{p}{q}$ (see Fig. 2).
Introduce the $a$-analogue of Fock function
\begin{equation}\label{focka}
\psi_a \left(\frac{p}{q}\right) {:=} \frac{1}{q} \arcosh \left(\frac{1}{2} X\left(\frac{p}{q}\right)\right).\\
\end{equation}

\medskip

\begin{theorem}
\label{main2}
The function $\psi_a$ can be extended to a continuous convex function of real $x \in [0,\frac{1}{2}]$, which is differentiable at all irrational $x$ and not differentiable at rational $x.$
\end{theorem}

We will explain now how to prove both  our theorems {(Theorem \ref{main} and Theorem \ref{main2})} by establishing links with known results about the stable norm and Mather $\beta$-function.

\medskip

\section{Stable norm and Mather $\beta$-function}

\subsection{The stable norm on $(\T^2,g)$}
In this section we introduce {Federer-Gromov's} {\it stable norm} on $H_1(\T^2;\R)$ and state some of its properties. We discuss only the case of the two-dimensional torus, which presents several simplifications; however, we refer the reader to \cite{BBI, GLP} for a more general presentation.

Let  $g$ be a Riemannian metric on $\T^2:=\R^2/\Z^2$ and let $L_g$ be the associated length functional.

For any integral class $h\in H_1(\T^2;\Z) \subset H_1(\T^2;\R)$ we define
$$
\ell_g(h) := \min\left\{L_g(\g):\; \g \; \mbox{is a smooth closed curve representing $h$}\right\}.
$$
Observe that $\ell_g(h)=\ell_g(-h)$, and $\ell_g(h)=0$ if and only if $h=0$.
It follows from a result by Hedlund \cite{Hedlund} (see also \cite[Theorem 8.5.10]{BBI}) that  $\ell_g: H_1(\T^2;\Z) \longrightarrow [0,+\infty)$ is positively homogeneous:
$$
\ell_g(n\,h) = n\,\ell_g(h) \qquad \forall \; h\in H_1(\T^2;\Z)\; {\rm and}\; \forall\; n\in \N;
$$
observe that this result is a peculiarity of the $2$-dimensional case.

The {\it stable norm} $\|\cdot\|_s$  corresponds to the unique norm on $H_1(\T^2;\R)$ which coincides with $\ell_g$ on $H_1(\T^2;\Z)$ (see \cite[Proposition 8.5.3]{BBI}). This norm can be actually constructed quite easily from $\ell_g$ in the following way: 
\begin{itemize}
\item[-] for every $h\in H_1(\T^2;\Z)$, we define $\|h\|_s:=\ell_g(h)$;
\item[-] then, using homogeneity, we extend it along lines of rational slopes: 
$$\|\alpha\, h\|_s := |\alpha| \,\|h\|_s \quad \forall\, h\in H_1(\T^2;\Z)
\; {\rm and}\; \forall\; \alpha\in \R;
$$
\item[-] using that $\ell_g(h_1 + h_2) \leq \ell_g(h_1) + \ell_g(h_2)$ for any $h_1, h_2\in H_1(\T^2;\Z)$,  then one can extend $\|\cdot\|_s$ continuously to the whole
$H_1(\T^2;\R)$.
\end{itemize}

%Here are some properties of $\|\cdot\|_s$.

The unit ball ${\mathcal B}_1:=\{h\in H_1(\T^2;\R): \, \|h\|_s\leq 1\}$ is strictly convex, namely its boundary does not contain straight line segments. Equivalently, if 
$h_1$ and $h_2$ are linearly independent, then
$$
\|h_1+h_2\|_s < \|h_1\|_s + \|h_1\|_s
$$
(see for example \cite[Exercise 8.5.15]{BBI}).

The following theorem was proven by Bangert \cite{Bangert2}  (see also related results by Mather \cite{Ma90} for twist maps, and by Klempnauer-Schr\"oder \cite{KS} for the case of Finsler metrics).

\begin{theorem}[Bangert]
Let $h=(h_1,h_2) \in H_1(\T^2;\R)\setminus \{0\}$, then
\begin{itemize}
\item[-] If $h_2\neq 0$ and $h_1/h_2 \in \R\setminus \Q$, then $\|\cdot\|_s$ is differentiable at $h$. 
\item[-] If $h_2=0$ or $h_2\neq 0$ and $h_1/h_2 \in  \Q$, then $\|\cdot\|_s$ is differentiable at $h$ if and only if
there exists a foliation of $\T^2$ by shortest closed geodesics in the same homotopy class, which is the primitive element in
$H_1(\T^2;\Z)\cap \{\R\cdot  h\}$.
\end{itemize}
\end{theorem}
%\item[] 

For some properties of minimal geodesics at irrational directions, see \cite{Schroeder}.
%\end{itemize}

\medskip

\subsection{Mather's $\beta$-function}
The stable norm is related to the so-called {\it Mather's minimal average action} (or {\it Mather's $\beta$ function}) \cite{Ma}. 
Hereafter we provide a brief introduction; we refer the reader to \cite{Sorrentinobook} for a more comprehensive one. 

Let us consider the Lagrangian associated  to the geodesic flow on $(\T^2,g)$:
\begin{eqnarray*}
L: T\T^2  &\longrightarrow& \R \\
(x,v) &\longmapsto& \frac{1}{2} g_x(v,v).
\end{eqnarray*}

Let ${\mathcal M}(L)$ the set of Borel probability measures $\mu$ on $T\T^2$ that are invariant under the Euler-Lagrange flow of $L$ ({\it i.e.}, the geodesic flow). To each element $\mu \in {\mathcal M}(L)$, we can associate its {\it homology} $\rho(\mu)\in H_1(\T^2;\R)$ (also called {\it Schwartzman asymptotic cycle}) in the following way  (see \cite{Ma, Sorrentinobook} for more details): it is the unique element of $H_1(\T^2;\R)$ for which 
\begin{equation}\label{homologymeasure}
\langle [\eta], \rho(\mu) \rangle = \int_{T\T^2} \eta(x,v)\, d\mu(x,v),
\end{equation}
where $\eta$ is any closed $1$-form on $\T^2$, $\eta(x,v)$ denotes the $1$-form $\eta$ computed at a tangent vector $(x,v)$, $[\eta]\in H^1(\T^2;\R)$ represents its cohomology, and $\langle \cdot, \cdot \rangle$ is the canonical pairing between (real) cohomology and homology groups. Observe that the integral on the right-hand side of \eqref{homologymeasure} only depends on the cohomology class of $\eta$ (see for example \cite[Lemma on p. 176]{Ma}).
One can prove that $\rho: {\mathcal M}(L) \longrightarrow H_1(\T^2;\R)$ is surjective and continuous \cite[Proposition 3.2.2]{Sorrentinobook}.

We define {\it Mather's $\beta$ function} as
\begin{eqnarray*}
\beta_{g}: H_1(\T^2;\R) &\longrightarrow& \R\\
h &\longmapsto& \beta_g(h):= \min_{\mu\in \rho^{-1}(h)} \int_{T\T^2} L_g(x,v) d\mu.
\end{eqnarray*}

There is the following relation with the stable norm  \cite[Proposition 1.4.2]{MassartThesis}:
$$\beta_g (h) =  \frac{1}{2} \|h\|_s^2.$$

\medskip

%%%%%%%%%%%%%%%%%%%%%%%%%%
\section{Fock's function and stable norm of punctured and one-holed tori}

We want to discuss how to prove  Theorem 2 by extending the above results {on Mather's $\beta$ function and the stable norm}, to  hyperbolic punctured and one-holed tori.\\ 

Let us consider first a punctured torus $\cT^2$, namely a surface homeomorphic to $\T^2\setminus\{{\rm point}\}$, equipped with a complete hyperbolic metric $g$ of finite area. The puncture corresponding to the removed point will be called {\it cusp} and hereafter denoted by  $C$.

The definition of the stable norm can be done as above. For each primitive homology class $h\in H_1(\cT^2;\Z)$, define $\ell_g(h)$ to be the length of the unique simple geodesic with homology class $h$ (see \cite[Corollary 1]{McShane}). For a non-primitive homology class $h$ there exists $n\in \N$ and $\hat{h}$ primitive such that  $h= n \hat{h}$ and we have
$\ell_g(h)= n \ell_g(\hat{h})$.
One can prove that $\ell_g$ satisfies the triangle inequality  (which is strict when we sum  linearly independent cycles).
The stable norm $\|\cdot\|_s$ can be then defined exactly as above: firstly, it coincides with $\ell_g$ on $H_1(\cT^2;\Z)$; then, it  is extended  homogeneously along lines of rational slopes, and  by continuity to the whole $H_1(\cT^2;\R)$.\\

The following theorem has been stated, without proof, in \cite{McShane} (we slightly rephrase its statement).

\begin{theorem}[McShane and Rivin]
Let $\|\cdot\|_s$ be the stable norm of $(\cT^2, g)$ and let $h=(h_1,h_2) \in H_1(\cT^2;\R)\setminus \{0\}$. Then the following hold true:
\begin{itemize}
\item[i)] If $h_2\neq 0$ and $h_1/h_2 \in \R\setminus \Q$, then $\|\cdot\|_s$ is differentiable at $h$ % and in particular is flat to infinite order. 
\item[ii)] If $h_2=0$ or $h_2\neq 0$ and $h_1/h_2 \in  \Q$, then $\|\cdot\|_s$ is not differentiable.
\end{itemize}
\end{theorem}
%\medskip

%\begin{Remark}{\rm
%Observe that differently from what happens on $\T^2$, the stable norm cannot be differentiable at rationals; in fact, there cannot be a foliation of $\cT^2$ consisting of (homologous) closed geodesics  (due to the hyperbolicity of the flow, there exists at most one simple closed geodesic in each homology class).}
%\end{Remark}

%\smallskip

Let us now describe how this result could be obtained from Bangert's results for $\T^2$. The idea is simply to ``plug the hole'' in a smooth way, which does not affect the minimizing properties of the metric. 
 
Recall that  any minimal compact lamination on a punctured torus does not intersect small cusp regions (namely, a neighbourhood of the cusp bounded by a horocycle); this result goes back at least to Poincar\'e, although a sharper version was proven by McShane \cite{McShaneThesis} (see also \cite[Theorems 1.1 \& 1.2]{McShaneRivin}).

\begin{prop}
Let $\eps>0$. Any punctured torus has a cusp region with bounding curve of length $4-\eps$ and this bound is optimal. No simple closed geodesic intersects a cusp region with boundary curve of length $4-\eps$.
\end{prop}

Therefore, we proceed as follows (cf. \cite{MPS}, proof of Corollary 1.7). First, remove from $\cT^2$ a cusp region of length, for example, equals to $1$; then, glue a euclidean hemisphere of equator length $1$. In this way, we obtain a two dimensional torus and the minimum length geodesics on this torus do not enter the added hemisphere. Therefore, the two metrics have the same stable norm. 

Applying now Bangert's theorem we have McShane-Rivin result. Indeed, there cannot be a foliation of $\T^2$ consisting of homologous closed geodesics since due to the hyperbolicity of the flow there exists at most one simple closed geodesic in each homology class.

Now let us consider a special case of hyperbolic equianharmonic punctured tori. In this case the lengths of simple closed geodesics are given by
$$
l=2\arcosh \frac{3}{2} m,
$$
where $m$  is the corresponding Markov number, and Fock's function $\psi(x)$ given by (\ref{fock}) is simply (half) the restriction of the corresponding stable norm to the line $q=1$ with $x=p/q$, where $(p,q)\in H_1(\cT^2; \mathbb Z)=\mathbb Z^2.$ By a general fact, the restriction of the norm to an affine line is a convex function on the line, so we have a proof of Fock's result (which is essentially the same as his own).
Now Theorem 1 follows from the above results of Bangert and McShane-Rivin.

The same arguments work for the hyperbolic tori with a hole, see \cite{MPS}. Indeed, we can cut the hyperbolic torus with a hole along the simple geodesic around the hole and glue the Euclidean hemisphere of the same equator length. It is  geometrically evident that any shortest simple closed geodesic on such  {``filled torus''} can not cross the hemisphere (cf. \cite{MPS}), so the stable norm remains the same. 

It is known after Fricke (see e.g. Goldman \cite{Goldman}) that the hyperbolic structures on one-holed torus are parametrized by  the positive real component of the cubic surface
$$
X^2+Y^2+Z^2-XYZ=c, \qquad c<0,
$$
where the geodesic length of the hole is given by $$l=2\arcosh \frac{2-c}{2}.$$
A natural action of the corresponding mapping class group $SL_2(\mathbb Z)$ is generated by the Vieta involutions
$(X, Y, Z) \; {\longmapsto} \; (X, Y, XY-Z)$ and cyclic permutations. This allows us to compute the corresponding stable norm recursively in the same way as for the Markov numbers.

Applying all this to the specific one-holed torus corresponding to the integer orbit (\ref{orbit}) we have the proof of Theorem 3.

\section{Concluding remarks} 

Let us mention some natural questions {that would need further investigation}.

 The picture of the unit ball for the corresponding stable norm from McShane-Rivin \cite{McShane} suggests that it has ``corners" at every rational points. Can we make this {more quantitative}? What can we say about the left and right derivatives of Fock's function at rational points, in particular at the symmetry point $x=\frac{1}{2}$? Since the function is convex they do exist. What are the corresponding values? Can we compute the derivatives of Fock's function at quadratic irrationals $x$? {We believe that}  Markov-Hurwitz ``most irrational" numbers (see e.g. \cite{Aigner,SV}) {might} play a special role {in this regard}.

%%%%%%%%%%%%%%%%%%%%%%%%%

\section{Acknowledgements} We are grateful to the University of Tel Aviv for the hospitality in the last week of April 2017, when the authors met and this work was essentially done, and to Misha Bialy and Leonid Polterovich for useful discussions. APV is very grateful also to Vladimir Fock for explanations of his results and to Leon Takhtajan for helpful discussions of Teichm\"uller spaces.

%%%%%%%%%%%%%%%%%%%%%%%%%

\end{document}